\documentclass[10pt]{article}
\usepackage{amsmath}
\usepackage{amssymb}
\usepackage{amsthm}
\usepackage{mathrsfs}
\usepackage{amsfonts}
\usepackage{tikz}
\usepackage{subfigure}

\newtheorem{theorem}{Theorem}[section]

\newtheorem{conjecture}[theorem]{Conjecture}

\title{A note on the Erd{\H o}s-Faber-Lov{\' a}sz Conjecture: quasigroups and complete digraphs\thanks{Research supported by: G. A-P. partially supported by CONACyT-M{\' e}xico under Projects 166306, 178395 and PAPIIT-M{\' e}xico under Project IN101912. C. R-M. partially supported by a CONACyT-M{\' e}xico Postdoctoral fellowship and by the National scholarship programme of the Slovak republic. A. V-{\' A}. partially supported by SNI of CONACyT-M{\' e}xico.}}

\author{Gabriela Araujo-Pardo \footnotemark[2]
\and Christian Rubio-Montiel \footnotemark[2] \footnotemark[3]
\and Adri{\' a}n V{\' a}zquez-{\' A}vila \footnotemark[4]
}

\begin{document}

\maketitle

\def\thefootnote{\fnsymbol{footnote}}
\footnotetext[2]{Instituto de Matem{\' a}ticas, Universidad Nacional Aut{\' o}noma de M{\' e}xico, 04510 M{\' e}xico City, Mexico, {\tt [garaujo|christian]@matem.unam.mx}.}
\footnotetext[3]{Department of Algebra, Comenius University, 842 48 Bratislava, Slovakia, {\tt christian.rubio@fmph.uniba.sk}.}
\footnotetext[4]{Subdirecci{\' o}n de Ingenier{\' i}a y Posgrado, Universidad Aeron{\' a}utica en Quer{\' e}taro, 76270 Quer{\' e}taro, Mexico, {\tt adrian.vazquez@unaq.edu.mx}.}


\begin{abstract}
A \emph{decomposition} of a simple graph $G$ is a pair $(G,P)$ where $P$ is a set of subgraphs of $G$, which partitions the edges of $G$ in the sense that every edge of $G$ belongs to exactly one subgraph in $P$. If the elements of $P$ are induced subgraphs then the decomposition is denoted by $[G,P]$.

A $k$-$P$-\emph{coloring} of a decomposition $(G,P)$ is a surjective function that assigns to the edges of $G$ a color from a $k$-set of colors, such that all edges of $H\in P$ have the same color, and, if $H_1,H_2\in P$ with $V(H_1)\cap V(H_2)\neq\emptyset$ then $E(H_1)$ and $E(H_2)$ have different colors. The \emph{chromatic index} $\chi'((G,P))$ of a decomposition $(G,P)$ is the smallest number $k$ for which there exists a $k$-$P$-coloring of $(G,P)$.

The well-known Erd{\H o}s-Faber-Lov{\'a}sz Conjecture states that any decomposition $[K_n,P]$ satisfies $\chi'([K_n,P])\leq n$. We use quasigroups and complete digraphs to give a new family of decompositions that satisfy the conjecture.
\end{abstract}




\section{Introduction}
Erd{\H o}s, Faber and Lov{\'a}sz, in 1972, conjectured the following (see \cite{MR0409246}): ``if $|A_i|=n$, $1\leq i \leq n$, and $|A_i\cap A_j|\leq1$, for $1\leq i < j \leq n$, then one can color the elements of the union $\bigcup_{i=1}^{n}A_i$ by $n$ colors, so that every set has elements of all the colors.'' This conjecture is called the Erd{\H o}s-Faber-Lov{\'a}sz Conjecture (for short EFL), and this can be set in terms of decompositions (see \cite{AV,MR2359282}).
\begin{conjecture}\label{conj:erdos:decompositions}
If $[K_{n},P]$ is a decomposition, then $\chi'([K_{n},P])\leq n$.
\end{conjecture}
In the following section we give a family of decompositions using finite quasigroups and complete digraphs satisfying Conjecture \ref{conj:erdos:decompositions}; this is a generalization of a previous result given in \cite{AV} and it is related with a result given in \cite{MR2359282}.

\section{Quasigroups and digraphs}
To begin with, we introduce definitions related to quasigroups, complete digraphs and linear-factorizations. A \emph{digraph} $D$ is a finite, non-empty set $V$ (the \emph{vertices} of $D$) together with a set $A$ of ordered pairs of elements of $V$ (the \emph{arcs} of $D$). We denote by $|V|$ the \emph{order} and by $|A|$ the \emph{size} of $D$ respectively.

A digraph $D$ is called \emph{symmetric} if whenever $(u,v)$ is an arc of $D$ then $(v,u)$ is an arc of $D$ --every graph can be interpreted as a symmetric digraph--. A \emph{directed cycle} or a \emph{$d$-gon} is a subdigraph with set of vertices $\{v_1,v_2,\dots,v_d\}$, such that their arcs are $(v_d,v_1)$ and $(v_i,v_{i+1})$ for $i\in\{1,\dots,d-1\}$ and $d\geq 2$. A \emph{loop} or a \emph{$1$-gon} is an arc joining a vertex with itself.

The \emph{complete digraph} $\overrightarrow{K}_n^*$ has order $n$ and size $n^2$ ($n$ loops and $\binom{n}{2}$ $2$-gons). A \emph{linear-factor} of the complete digraph $\overrightarrow{K}_n^*$ is a subdigraph of order $n$ and size $n$, such that it is a set of pairwise vertex-disjoint $d$-gons. A \emph{linear-factorization} of $\overrightarrow{K}_n^*$ is a set of pairwise arc-disjoint linear-factors, such that these linear-factors induce a partition of the arcs, see Figure \ref{Fig1}: c).

A \emph{quasigroup $(\mathcal{Q}_n,\cdot)$} is a set $\mathcal{Q}$ of $n$ elements with a binary operation $\cdot$, such that for each $x$ and $y$ in $\mathcal{Q}$ there exist unique elements $a$ and $b$ in $\mathcal{Q}$ with $x\cdot a=y$ and $b\cdot x=y$.

Let $(\mathcal{Q}_n,\cdot )$ be a quasigroup and the complete digraph $\overrightarrow{K}_n^*$, such that its vertices are the elements of $\mathcal{Q}_n$. Afterwards, we color the arcs of $\overrightarrow{K}_n^*$ by $n$ colors which are in a one-to-one correspondence with the elements of $\mathcal{Q}_n$ so that for any two vertices $x$ and $y$ in $\mathcal{Q}_n$ the arc $(x,y)$ obtains the color corresponding to $a\in\mathcal{Q}_n$ for which $x\cdot a=y$ holds true. Then the resulting graph with the described coloring of arcs is called the \emph{Cayley color graph} $C(\mathcal{Q}_n)$ of $\mathcal{Q}_n$. The Cayley color graph of a quasigroup is described in \cite{MR0360349}.

It is not hard to prove that the arcs colored by the same color in $C(\mathcal{Q}_n)$ induce a linear-factor of this digraph. An arc colored by the color corresponding color to some $a\in \mathcal{Q}_n$ outgoing from the vertex $x$ leads into $x\cdot a$ in $C(\mathcal{Q}_n)$. The element $x\cdot a$ is exactly one for any $x$ and any $a$ of $\mathcal{Q}_n$.

Consequently, the Cayley color graph $C(\mathcal{Q}_n)$ can be considered as a linear-facto\-rization $\mathcal{F}$ of $\overrightarrow{K}_n^*$ of $n$ linear-factors. In \cite{MR0360349} it was proved that any linear-factorization $\mathcal{F}$ of the complete digraph $\overrightarrow{K}_n^*$ and any one-to-one mapping of the vertex set of $\overrightarrow{K}_n^*$ onto the set of linear-factors of $\mathcal{F}$ determines a quasigroup $\mathcal{Q}_n$, such that the Cayley color graph $C(\mathcal{Q}_n)$ of $\mathcal{Q}_n$ can be considered $(\overrightarrow{K}_n^*,\mathcal{F})$, as described above.

Following, we relate the previous concepts with decompositions of complete graphs. Let $[K_n,P]$ be a decomposition $P$ of $K_n$ and let $\overrightarrow{K}_n$ be the symmetric complete digraph (without loops). We consider the decomposition $[\overrightarrow{K}_n,P]$ induced by $[K_n,P]$, that is, $P$ is a set of subdigraphs of $\overrightarrow{K}_n$, which partitions the arcs of $\overrightarrow{K}_n$ in the sense that every arc of $\overrightarrow{K}_n$ belongs to exactly one subdigraph in $P$ and every element of $P$ is a symmetric complete subdigraph. The digraph $\overrightarrow{K}_n^*$ is $\overrightarrow{K}_n$ with the set $L$ of $n$ loops.

Now, we state and prove the main theorem:
\begin{theorem}\label{teoAssgnment}
Let $[\overrightarrow{K}_n,P]$ be a decomposition $P$ of $\overrightarrow{K}_n$ arising from $[K_{n},P]$ and let $(\overrightarrow{K}_n^*,\mathcal{F})$ be a linear-factorization $\mathcal{F}$ of $\overrightarrow{K}_n^*$. If there exists a function $h\colon P \rightarrow \mathcal{F}$, such that for any $p\in P$, $(A(p)\cup L) \cap A(h(p))$ is a linear-factor $F_p$ of $p^*$ --$p$ with loops-- and for any $p,q\in P$, $A(F_p)\cap A(F_q)=\emptyset$ then $\chi'([K_{n},P])\leq n$.
\begin{proof}
Color the edges of an element $p$ of $P$ with $f(h(p))$ where $f$ is a one-to-one mapping of a quasigroup $\mathcal{Q}$ onto the set of linear-factors of $\mathcal{F}$. The $n$-coloring is well-defined due to the fact that for any $p,q\in P$, $A(F_p)\cap A(F_q)=\emptyset$ and the result follows.
\end{proof}
\end{theorem}
We can explain Theorem \ref{teoAssgnment} as following: 

\begin{figure}
  \begin{center}
    \includegraphics{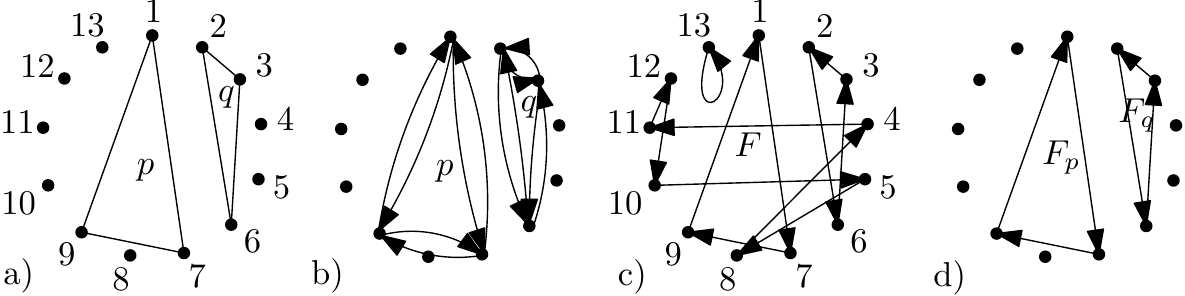}
    \caption{a) Two elements $p$ and $q$ of a decomposition of $K_{13}$ into triangle arising from the cyclic Steiner System $STS(13)$. b) $K_{13}$ as a symmetric digraph c) A linear-factor $F$ for $n=13$. The mapping $i\mapsto i+1$ produces a linear-factorization. d) The restriction of $F$ onto $p$ and $q$.}
    \label{Fig1}
  \end{center}
\end{figure}

Let $(\overrightarrow{K}_n^*,\mathcal{F})$ be a linear-factorization $\mathcal{F}$ of $\overrightarrow{K}_n^*$. Then every decomposition $P$ formed by complete subdigraphs obtained via some linear-factor $f_0$ of $\mathcal{F}$, meaning, the intersection of the arcs of $p\in P$ with the arcs of $f_0$ is a linear factor of $p$ has a consequence that $\chi'([K_{n},P])\leq n$. Figure \ref{Fig1} illustrates Theorem \ref{teoAssgnment} with an example for $n=13$. 

\bibliographystyle{amsplain}
\bibliography{biblio}
\end{document}